\documentclass[11pt]{article}
\usepackage{amsmath}
 \usepackage{amssymb}
 \usepackage{array}
 \usepackage{geometry}
 \usepackage{graphicx}
 \usepackage{enumerate}
\newtheorem {thm}{Theorem}[section]
\newtheorem{lem}{Lemma}[section]

\newtheorem {ex}{Example}[section]
\newtheorem {de}{Definition}[section]
\numberwithin{equation}{section} 
\begin{document}
\begin{center}
{\Large \textbf{On the Levy density function}}\\
\end{center}

{\begin{center} Jung Hun Han   \footnote{Corresponding Author Address: Centre for Mathematical Sciences, Pala Campus,\\
Arunapuram P.O., Pala, Kerala-686 574, India, ~~Email :
jhan176@gmail.com}\end{center}}

\begin{abstract}
In this paper, we introduce the Levy density function as the limit
of a generalized Mittag-Leffler density function. The fractional
integral equation for the generalized Mittag-Leffler density
function is also given. And the role of the Levy structure in the
fractional calculus is described. Finally, a transformation is
defined.
\end{abstract}

keywords: density function, fractional calculus, Mellin
convolution operator\\
$~~~~~~~~~~~~~~~~~~~~~~$written April, 30, 2010

\section{\bf  Introduction}\label{s1}
In the following sections, we are going to look into some details
related to Levy distribution and its role in the fractional
calculus. We state some preliminary results. Pochammer symbol is
defined as
\begin{equation}\label{} (b)_{k}=b(b+1)\cdots (b+k-1), ~(b)_{0}=1,~ b \neq 0. \end{equation}

\begin{de}
The  Euler gamma function $\Gamma (z)$ is defined as follows
\begin{eqnarray}
  \Gamma (z) &=& p^{z}\int^{\infty}_{0}t^{z-1}e^{-pt}dt, ~~\mathfrak{R}(p)>0, \mathfrak{R}(z)>0 \\
   &=&\lim_{n\rightarrow\infty} \frac{n! n^{z}}{z(z+1)\cdots (z+n)},~~z\neq 0, 1, 2, 3, \cdots.
   \end{eqnarray}
\end{de}
The representation of the Pochammer symbol in terms of gamma
functions is
\begin{equation}\label{} (a)_{k}= \frac{\Gamma(a+k)}{\Gamma(a)} \end{equation}
whenever the gammas exist.
\begin{lem}\label{hhh15}
[ Stirling asymptotic formula ] \cite{ma3, mh3}\\
For $|z|\rightarrow \infty$ and $\alpha$ a bounded quantity,
\begin{equation}\label{} \Gamma (z+ \alpha) \approx (2\pi)^{1/2}z^{z+\alpha
-1/2}e^{-z}.
\end{equation}
\end{lem}

\begin{lem}\label{hhh15}

For $|\gamma|\rightarrow \infty$ and $k$ a bounded quantity,
\begin{equation}\label{} \lim_{\gamma\rightarrow\infty}\frac{(\gamma)_{k}}{\gamma^{k}}=
\lim_{\gamma\rightarrow\infty}
\frac{\Gamma(\gamma+k)}{\Gamma(\gamma)\gamma^{k}} =1.
\end{equation}
\end{lem}

$H$-function representations and their convergent regions are
important. We are not going to explain them and the definition of
Mittag-Leffler function, its generalized ones and their properties
are not added here. So for those, \cite{ma3, mh3, mhs2} have
detailed descriptions.

In the paper \cite{ma1}, the author used some properties of Mellin
transformation and statistical methods to investigate the
Mittag-Leffler statistical distribution and its generalized ones.
In \cite{ma1}, he showed the possibility of getting the explicit
form of Levy density function. He also computed
$E[v^{\frac{1}{\alpha}}]^{s-1}$ from a gamma density (\cite{ma1}
for details). This computation will lead us to define a
transformation in section 4. In \cite{mh2}, the following results
are stated.

\begin{lem} For $\Re(\beta)>0,\Re(\gamma)>0$,
 \begin{equation}\label{}
 \lim_{|\beta|\rightarrow \infty} \frac{\Gamma(\beta)}{\Gamma(\gamma)}
 H^{1,1}_{1,2}\left[-(z\beta^{\frac{\alpha}{\delta}})^{\delta}
 \big|^{(1-\gamma,1)}_{(0,1),(1-\beta,\alpha)}\right]
 =[1-z^{\delta}]^{-\gamma}
\end{equation}
\end{lem}

\begin{thm}\label{t1} For $\Re(\beta)>0,\Re(\gamma)>0,x>0,a>0,q>1,c>0$,
 \begin{equation}\label{}
 \lim_{|\beta|\rightarrow \infty} \frac{\Gamma(\beta)}{\Gamma(\frac{\eta}{q-1})}
 x^{\gamma}E^{\frac{\eta}{q-1}}_{(\alpha,\beta)}\left[-a(q-1)(x\beta^{\frac{\alpha}{\delta}})^{\delta}
\right]=cx^{\gamma}\left[1+a(q-1)x^{\delta}
\right]^{-\frac{\eta}{q-1}}.
\end{equation}
\end{thm}
Theorem \ref{t1} is related to the pathway model, for $x>0$ from
where one has Tsallis statistics, superstatistics, power law and
many others.

The concept of pathway model is important in the sense that it is
a pathway between two totally different looking systems. It can be
considered that the behavior of the family of density functions
tied up by the pathway parameter is similar to that of the
collection of germs connecting two different points via a path in
the complex plane. For details, see \cite{ga, ma1, mh1, mh2}.
\\
\\
In this paper, we construct the Levy density function through a
pathway parameter. Finally we define a transformation of
connecting ordinary space and $\alpha$-fractional (or
$\alpha$-level) space.


\section{\bf The Levy density function}\label{s3}

Consider a gamma random variable $x$ with the density function
\begin{equation}\label{hhh18}
g(x)=\begin{cases} \frac{\gamma^{\gamma}x^{\gamma-1}e^{-\gamma
x}}{\Gamma(\gamma)}~~\text{for}~
0< x<\infty , ~\gamma >0\\
0, ~~~\text{otherwise.} \end{cases}\end{equation} This density
function has some interesting properties.
\begin{enumerate}
    \item Laplace transform : $L_{g}(s)=\left(1+\frac{s}{\gamma}\right)^{-\gamma}$
    \item Mellin transform : $M_{g}(s)=\frac{\Gamma(\gamma+s-1)
}{\Gamma(\gamma)} \gamma^{1-s}$
    \item Mellin-Barnes integral representation : $$\frac{1}{2\pi i}\oint_{L'}
\frac{\Gamma(\gamma+s-1) }{\Gamma(\gamma)} \gamma^{1-s}x^{-s}ds$$
\end{enumerate}

 We consider
\begin{equation}
\label{hhh13}f(x)=x^{\alpha\gamma-1}\gamma^{\gamma}E^{\gamma}_{(\alpha,\alpha\gamma)}
(-\gamma x^{\alpha}),\end{equation} which can be obtained by using
the technique developed in section 4 with the density function in
\eqref{hhh18}. This is a generalized Mittag-Leffler density
function and the function which leads to the explicit form of Levy
density function.

\begin{enumerate}
    \item Laplace transform : \begin{equation}\label{hhh12}L_{f}(s)=\left(1+\frac{s^{\alpha}}
    {\gamma}\right)^{-\gamma}.\end{equation}
    \item Mellin transform : \begin{equation}\label{hhh22}
M_{f}(s)=\frac{\Gamma(\gamma+\frac{s}{\alpha}-\frac{1}{\alpha})
\Gamma(-\frac{s}{\alpha}+\frac{1}{\alpha})}{\alpha\Gamma(\gamma)\Gamma(1-s)}
\gamma^{-\frac{s}{\alpha}+\frac{1}{\alpha}}\end{equation} where
$0< Re(s)<\alpha<1$.
    \item Mellin-Barnes integral representation : \begin{equation}\frac{1}{2\pi i}\oint_{L}
\frac{\Gamma(s)\Gamma(\gamma-s)}{\Gamma(\gamma)\Gamma(\alpha\gamma-\alpha
s)} x^{\alpha\gamma-\alpha s-1}\gamma^{-s+\gamma}ds\end{equation}
\begin{equation}\label{hhh14}
=\frac{1}{2\pi i}\oint_{L'}
\frac{\Gamma(\gamma+\frac{s}{\alpha}-\frac{1}{\alpha})
\Gamma(-\frac{s}{\alpha}+\frac{1}{\alpha})}{\alpha\Gamma(\gamma)\Gamma(1-s)}
\gamma^{-\frac{s}{\alpha}+\frac{1}{\alpha}}x^{-s}ds\end{equation}
where $0< Re(s)<\alpha<1$.
\end{enumerate}

Firstly, we look into its Laplace transform. When $\gamma$ tends
to $\infty$ in \eqref{hhh12},  \eqref{hhh12} goes to
$e^{-t^{\alpha}}$ which is the Laplace transform of Levy
distribution. This means that \eqref{hhh13} will lead us to the
Levy density function. For this, we make use of \eqref{hhh14} and
lemma \ref{hhh15}. From the representation \eqref{hhh14}, we have
$$
\frac{1}{2\pi i}\oint_{L'}
\frac{\Gamma(\gamma+\frac{s}{\alpha}-\frac{1}{\alpha})
\Gamma(-\frac{s}{\alpha}+\frac{1}{\alpha})}{\alpha\Gamma(\gamma)\Gamma(1-s)}
\gamma^{-\frac{s}{\alpha}+\frac{1}{\alpha}}x^{-s}ds.$$ As $\gamma$
tends to $\infty$, we get
\begin{equation}\label{}
\frac{1}{2\pi i}\oint_{L'}
\frac{\Gamma(-\frac{s}{\alpha}+\frac{1}{\alpha})}{\alpha\Gamma(1-s)}
x^{-s}ds, ~ 0< Re(s)<\alpha<1 \end{equation} which is the
Mellin-Barnes integral representation of Levy density function
\cite{ma1}.

\begin{de}
Let f(x) have $H$-function representation and be convergent
\cite{ma3, mhs2}. Then f(x) is said to be a function with Levy
structure if its $H$-function representation has the factor
$\frac{\Gamma(-\frac{s}{\alpha}+\frac{1}{\alpha})}{\alpha\Gamma(1-s)}$
in the integrand. In short,
$\frac{\Gamma(-\frac{s}{\alpha}+\frac{1}{\alpha})}{\alpha\Gamma(1-s)}$
can be called Levy structure.
\end{de}

In \cite{ma4}, the 2- parameter Weibull density function is
defined as follows,
\begin{equation}\label{}
f_1(x)=\begin{cases} \delta b x^{\delta-1}
e^{-bx^{\delta}}~~\text{for}~
0\leq x<\infty , ~\delta >0, ~b>0 \\
0, ~~~\text{elsewhere.} \end{cases}.\end{equation}

Consider the remaining part of \eqref{hhh14} after removing the
Levy structure,
\begin{equation}\label{han14}
\frac{1}{2\pi i}\oint_{L'}
\frac{\Gamma(\gamma+\frac{s}{\alpha}-\frac{1}{\alpha})
}{\Gamma(\gamma)} \gamma^{-\frac{s}{\alpha}+\frac{1}{\alpha}}
x^{-s}ds.\end{equation} Then by the residue theorem, it becomes
 $\frac{\alpha \gamma^{\gamma}}{\Gamma(\gamma)}
x^{\alpha\gamma-1}e^{-\gamma x^{\alpha}}$, where $\alpha >0,
\gamma
> 0$, which looks the same
as $\frac{\gamma^{\gamma}x^{\gamma-1}e^{-\gamma
x}}{\Gamma(\gamma)}$ when $\alpha$ is replaced by $1$. This is a
density function  and a  2-parameter generalized gamma density
function. we just list 3 cases in the below table. \vskip 0.2cm
\noindent\emph{\textbf{generalized gamma density functions from
Mittag-Leffler
density functions}}\\
\begin{tabular}{|c|c|c|c|}
  \hline
   generalized gamma density functions  &Mittag-Leffler
density functions    \\
  \hline
    $\alpha
x^{\alpha-1}e^{- x^{\alpha}}$ & $x^{\alpha-1}E_{(\alpha,\alpha)} (- x^{\alpha}) $    \\
\hline
    $\frac{\alpha \gamma^{\gamma}}{\Gamma(\gamma)}
x^{\alpha\gamma-1}e^{-\gamma x^{\alpha}}$ &
$x^{\alpha\gamma-1}\gamma^{\gamma}E^{\gamma}_{(\alpha,\alpha\gamma)}
(-\gamma x^{\alpha})$   \\
\hline
    $\frac{\alpha \delta^{\eta}}{\Gamma(\eta)}
x^{\alpha\eta-1}e^{-\delta x^{\alpha}}$ &
$x^{\alpha\eta-1}\delta^{\eta}E^{\eta}_{(\alpha,\alpha\eta)}
(-\delta x^{\alpha})$   \\
  \hline
\end{tabular}
\\

 When
$\gamma\rightarrow \infty$ in \eqref{han14}, the Mellin transform
becomes $1$, which means
\begin{equation}\label{}
\lim_{\gamma\rightarrow\infty}
\int^{\infty}_{0}x^{s-1}\frac{\alpha
\gamma^{\gamma}}{\Gamma(\gamma)} x^{\alpha\gamma-1}e^{-\gamma
x^{\alpha}}dx=
\lim_{\gamma\rightarrow\infty}\frac{\Gamma(\gamma+\frac{s}{\alpha}-\frac{1}{\alpha})
}{\Gamma(\gamma)}
\gamma^{-\frac{s}{\alpha}+\frac{1}{\alpha}}=1.\end{equation} So as
$\gamma$ takes very large value, then the parameter $s$ becomes
redundant.


\section{\bf  The fractional integral equation as an extension of the reaction rate model }\label{s4}
In \cite{mh2}, they consider fractional integral equations as
extensions of the reaction rate model.
$\frac{dN(x)}{dx}=-cN(x),~c>0 \Rightarrow N(x)-N_{0}= - c \int
N(x)dx $~ : reaction rate model which gives
     $N(x)-N_{0}= - c^{\alpha}~ _{0}D^{-\alpha}_{x}N(x)$
    and from there ~~
\begin{equation}
  \label{hhh1} N(x)-N_{0}f(x)= - c^{\alpha}~
  _{0}D^{-\alpha}_{x}N(x)
  \end{equation}
where $f(x)$ is a general integrable function on the finite
interval $[0,b]$.

By applying Laplace transformation,
$$ \widetilde{N(s)}-N_{0}
F(s)= - c^{\alpha}~ s^{\alpha} \widetilde{N(s)}$$ where $F(s)$ is
the Laplace transform of $f(x)$ and $\widetilde{N(s)}$ is the
Laplace transform of $N(x)$. Then

$$\widetilde{N(s)}=\frac{N_{0}s^{\alpha}
F(s)}{s^{\alpha}+c^{\alpha}}=N_{0} F(s)\sum^{\infty}_{k=0}(-1)^{k}
c^{\alpha k} s^{-\alpha k}$$
$$=N_{0} F(s)+N_{0} F(s)\sum^{\infty}_{k=1}(-1)^{k} c^{\alpha k}
s^{-\alpha k}=N_{0} F(s)-N_{0} F(s)\sum^{\infty}_{k=0}(-1)^{k}
c^{\alpha k+\alpha} s^{-\alpha k -\alpha}.$$ But note that
$\sum^{\infty}_{k=0}(-1)^{k} c^{\alpha k+\alpha} s^{-\alpha k
-\alpha}$ is the Laplace transform of
$x^{\alpha-1}E_{(\alpha,\alpha)}(-c^{\alpha}x^{\alpha})$. By
applying the inverse Laplace transformation,
\begin{eqnarray}\label{hhh10}
N(x)&=&N_{0} f(x) +
N_{0}\sum^{\infty}_{k=1}\frac{(-1)^{k}c^{\alpha
k}}{\Gamma(\alpha k)} \int^{x}_{0}(x-t)^{\alpha k-1}f(t)dt\\
\label{hhh20} &=&N_{0} f(x)-N_{0}c^{\alpha}
\int^{x}_{0}(x-t)^{\alpha-1}E_{(\alpha,\alpha)}(-c^{\alpha}(x-t)^{\alpha})f(t)dt.
\end{eqnarray}

It seems that the kernel of the integral in \eqref{hhh20} has Levy
structure and  that if the kernel of an integral operator has Levy
structure, it goes well with the function possessing Levy
structure, in other words, it brings more nice forms and if $f(x)$
has the Levy structure, then the solution has the Levy structure,
too. Here we give some examples.\\
Let $f(x)$ be $1$ in \eqref{hhh1}, then
$N(x)=N_{0}E_{\alpha}(-c^{\alpha}x^{\alpha})$ from \eqref{hhh20}.
This can be written as follows. \vskip 0.1cm
\textbf{\emph{functions without Levy structure}}\\
\begin{tabular}{|c|c|c|}
  \hline
  $f(x)$ &$\rightarrow$  & $N(x)$ \\
  \hline
  $1$ &  & $N_{0}E_{\alpha}(-c^{\alpha}x^{\alpha})$ \\
   \hline
   $x$&  &$N_{0}xE_{(\alpha,2)}(-c^{\alpha}x^{\alpha})$  \\
   \hline
  $e^{-x}$ &  & $N_{0}\sum^{\infty}_{k=0}(-x)^{k}E_{(\alpha,k+1)}(-c^{\alpha}x^{\alpha})$ \\
   \hline
   $ e^{-(cx)^{\alpha}}$&  &$N_{0}\sum^{\infty}_{k=0}(-1)^{k}x^{k\alpha} c^{k\alpha}
\frac{\Gamma(\alpha k +1)}{k!}
E_{(\alpha,k\alpha+1)}(-c^{\alpha}x^{\alpha})$  \\
  \hline
   $\frac{x^{\mu-1}}{\Gamma(\mu)}$ &  & $N_{0}x^{\mu-1}E_{(\alpha,\mu)}(-c^{\alpha}x^{\alpha})$ \\
   \hline
   $x^{\mu-1}E^{\gamma}_{(\alpha,\mu)}(-c^{\alpha}x^{\alpha})$&  &$N_{0}x^{\mu-1}E^{\gamma+1}_{(\alpha,\mu)}(-c^{\alpha}x^{\alpha})$  \\
   \hline
\end{tabular}
\\

\vskip 0.2cm
\textbf{\emph{functions with Levy structure}}\\
\begin{tabular}{|c|c|c|}
  \hline
  $f(x)$ &$\rightarrow$  & $N(x)$ \\
   \hline
   $\frac{x^{\alpha-1}}{\Gamma(\alpha)}$ &  & $N_{0}x^{\alpha-1}E_{(\alpha,\alpha)}(-c^{\alpha}x^{\alpha})$ \\
   \hline
   $x^{\alpha-1}E^{\gamma}_{(\alpha,\alpha)}(-c^{\alpha}x^{\alpha})$&  &$N_{0}x^{\alpha-1}
   E^{\gamma+1}_{(\alpha,\alpha)}(-c^{\alpha}x^{\alpha})$  \\
   \hline
  \end{tabular}
  \vskip 0.1cm

Now, we are in a position to show that the density function
\eqref{hhh13} has its own integral equation. By changing
$c^{\alpha}$ to $\gamma$,  \eqref{hhh1} becomes
\begin{equation}\label{}
N(x)-N_{0}f(x)= - \gamma ~_{0}D^{-\alpha}_{x}N(x).\end{equation}
Then from \eqref{hhh20}, the solution is accordingly
\begin{eqnarray*}\label{}
N(x)&=&N_{0} f(x)-N_{0}\gamma
\int^{x}_{0}(x-t)^{\alpha-1}E_{(\alpha,\alpha)}(-\gamma(x-t)^{\alpha})f(t)dt \\
&=&N_{0}f(x)- N_{0}\gamma \int^{x}_{0}
(x-t)^{\alpha-1}H^{1,1}_{1,2}\left[
\gamma(x-t)^{\alpha}\big|^{(0,1)}_{(0,1)(1-\alpha,\alpha)}
   \right] f(t)dt \\
&=& N_{0} f(x)-N_{0} \int^{x}_{0}\frac{1}{2\pi i} \oint_{L}
\frac{\Gamma(1+\frac{s}{\alpha}
-\frac{1}{\alpha})\Gamma(\frac{1}{\alpha}- \frac{s}{\alpha})
}{\alpha \Gamma(1-s)}(x-t)^{-s}\gamma^{\frac{1-s}{\alpha}}ds f(t)
dt.\end{eqnarray*}

If we consider this particular equation, then $
x^{\alpha\gamma-1}\gamma^{\gamma}E^{\gamma}_{(\alpha,\alpha\gamma)}
(-\gamma x^{\alpha})$ becomes the solution of
\begin{equation}\label{}
N(x)-N_{0}x^{\alpha\gamma-1}\gamma^{\gamma}E^{\gamma-1}_{(\alpha,\alpha\gamma)}
(-\gamma x^{\alpha})= - \gamma
~_{0}D^{-\alpha}_{x}N(x),\end{equation} which is the main
generalized Mittag-Leffler density function in section 2.

\section{\bf A lifting from ordinary space to $\alpha$-fractional($\alpha$-level) space}\label{s5}

In \cite{ma1}, he shows a process to lift a gamma density to a
generalized Mittag-Leffler density function by using statistical
techniques and it is in Example \ref{hhh31}.
\begin{ex}\label{hhh31}
Let $x$ be a simple exponential random variable with the density
function $f(x)=e^{-x}$. We attach the Levy structure to $
E[x^{\frac{1}{\alpha}}]^{s-1} =\int^{\infty}_{0}
(x^{\frac{1}{\alpha}})^{s-1}e^{-x}dx$.

  Then
\begin{equation}
   x^{\alpha-1}E_{(\alpha,\alpha)} (- x^{\alpha})=\frac{1}{2\pi i} \oint_{L}
\frac{E[x^{\frac{1}{\alpha}}]^{s-1}\Gamma(-\frac{s}{\alpha}+\frac{1}{\alpha})x^{-s}}
{\alpha\Gamma(1-s)}ds
  \end{equation}
by the residue theorem.
\end{ex}
In this section, we will show through a new process how to lift a
function in the ordinary space to a corresponding function in the
$\alpha$-fractional(or $\alpha$-level) space.

 Let $f(x)$ be a function, which does not have Levy structure and lives in the
ordinary space and $h(x)$ be the Levy density function as a
kernel.
\begin{de}
 Define a transform of $f(x)$

\begin{eqnarray}
 J(f)(x) = \lim_{\gamma\rightarrow\infty}\int^{x}_{0} \left(\frac{t}{x}\right)f_{2}\left(\left(\frac{x}{t}\right)^{\alpha}\right)
 \frac{t^{\alpha\gamma-1}\gamma^{\gamma}E^{\gamma}_{(\alpha,\alpha\gamma)}
(-\gamma t^{\alpha})}{t}dt
   \end{eqnarray}
   where ~$\alpha$ fixed in $0 < \alpha < 1$, x>0, $f_{2}(x)=xf(x)$ and $f$ is integrable and continuous on
   the interval.
\end{de}

\begin{ex}\label{hhh16}
Let $x$ be a exponential random variable with the density function
\begin{equation}\label{}
\nonumber f(x)=e^{-x},~
f_{2}(x)=x\sum^{\infty}_{l=0}\frac{(-1)^{l}x^{l}}{l!}.\end{equation}
Then
\begin{eqnarray*}
  J(f)(x) &=&\lim_{\gamma\rightarrow\infty}\int^{x}_{0} \left(\frac{t}{x}\right)f_{2}\left(\left(\frac{x}{t}\right)^{\alpha}\right)
 \frac{t^{\alpha\gamma-1}\gamma^{\gamma}E^{\gamma}_{(\alpha,\alpha\gamma)}
(-\gamma t^{\alpha})}{t}dt\\
   &=& \lim_{\gamma\rightarrow\infty}\int^{x}_{0}
   \left(\frac{x}{t}\right)^{\alpha -1}\sum^{\infty}_{l=0}
   \frac{(-1)^{l}\left(\frac{x}{t}\right)^{\alpha l}}{l!}
t^{\alpha\gamma-1}\gamma^{\gamma}E^{\gamma}_{(\alpha,\alpha\gamma)}
(-\gamma t^{\alpha}) \frac{dt}{t} \\
 &=&\lim_{\gamma\rightarrow\infty}\int^{1}_{0} u^{1-\alpha}
\sum^{\infty}_{l=0} \frac{(-1)^{l}u^{-\alpha l}}{l!}
   \left(x u\right)^{\alpha\gamma-1}\gamma^{\gamma}\sum^{\infty}_{k=0}\frac{(\gamma)_{k}
(-1)^{k}\gamma^{k} \left(x u\right)^{\alpha k}}{k!\Gamma(\alpha
k+\alpha\gamma)}
 \frac{du}{u}\\
  \end{eqnarray*}
\begin{eqnarray*}
&=&\lim_{\gamma\rightarrow\infty}\sum^{\infty}_{l=0}
   \frac{(-1)^{l}}{l!}\gamma^{\gamma}
\sum^{\infty}_{k=0}\frac{(\gamma)_{k} (-1)^{k}\gamma^{k} x^{\alpha
k +\alpha\gamma-1} }{k!\Gamma(\alpha k+\alpha\gamma)}
 \int^{1}_{0}u^{\alpha\gamma +\alpha k-\alpha l-\alpha-1}du\\
 &&(\mbox{fix $\gamma$ so that the integral exists,})\\
&=&\lim_{\gamma\rightarrow\infty} \sum^{\infty}_{l=0}
   \frac{(-1)^{l}}{l!} \gamma^{\gamma}
\sum^{\infty}_{k=0}\frac{(\gamma)_{k} (-1)^{k}\gamma^{k} x^{\alpha
k +\alpha\gamma-1} }{k!\Gamma(\alpha k+\alpha\gamma)}
 \frac{\Gamma(\alpha k+\alpha\gamma-\alpha l -\alpha)}
 {\Gamma(\alpha k+\alpha\gamma-\alpha l-\alpha+1)}\\
 &=&\lim_{\gamma\rightarrow\infty}\sum^{\infty}_{l=0}
   \frac{(-1)^{l}}{l!} \gamma^{\gamma} \frac{1}{2\pi i} \oint_{L}
\frac{\Gamma(s)\Gamma(\gamma-s)\Gamma(-\alpha
s+\alpha\gamma-\alpha l -\alpha)\gamma^{-s}x^{-\alpha s
+\alpha\gamma-1}}{\Gamma(\gamma)\Gamma(\alpha\gamma-\alpha
s)\Gamma(-\alpha
s+\alpha\gamma-\alpha l -\alpha+1)}ds\\
&=&\lim_{\gamma\rightarrow\infty}\sum^{\infty}_{l=0}
   \frac{(-1)^{l}}{l!}\frac{1}{2\pi i}
\oint_{L}
\frac{\Gamma(\gamma+\frac{s}{\alpha}-\frac{1}{\alpha})\Gamma(-\frac{s}{\alpha}+\frac{1}{\alpha})\Gamma(1-
s-\alpha
l-\alpha)\gamma^{-\frac{s}{\alpha}+\frac{1}{\alpha}}x^{-s}}
{\Gamma(\gamma)\alpha\Gamma(1-s)\Gamma(2- s-\alpha
l -\alpha)}ds\\
&&(\mbox{where } 0 < R(s) < \alpha < 1)\\
&=&\sum^{\infty}_{l=0}
   \frac{(-1)^{l}}{l!} \frac{1}{2\pi i} \oint_{L}
\frac{\Gamma(-\frac{s}{\alpha}+\frac{1}{\alpha})\Gamma(1- s-\alpha
l -\alpha)x^{-s}} {\alpha\Gamma(1-s)\Gamma(2- s-\alpha l-\alpha)}ds\\
&&\frac{\Gamma(-\frac{s}{\alpha}+\frac{1}{\alpha})}{\alpha\Gamma(1-s)}\mbox{
has no poles in the interval }  0 < R(s) < \alpha < 1,\\
&&\mbox{ so we proceed with } \frac{\Gamma(1- s-\alpha
l-\alpha)}{\Gamma(2- s-\alpha l
-\alpha)} \mbox{by putting } 1-s-\alpha l-\alpha=-\alpha s_{1},\\
&=&\sum^{\infty}_{l=0}
   \frac{(-1)^{l}}{l!} \frac{1}{2\pi i} \oint_{L}
\frac{\Gamma(l+1-\alpha s_{1})\Gamma(-\alpha s_{1})x^{-\alpha
s_{1} +\alpha
l+\alpha-1}} {\Gamma(\alpha l+\alpha-\alpha s_{1})\Gamma(1-\alpha s_{1})}ds\\
&&\mbox{the integrand can have a pole at 0 and only one pole there
since }0<\alpha<1 \\
&=&\sum^{\infty}_{l=0}
   \frac{(-1)^{l}}{l!}  \frac{\Gamma(1+l)x^{\alpha-1+\alpha l}}
{\Gamma(\alpha l+\alpha)}\\
&=&x^{\alpha-1}E_{(\alpha,\alpha)} (- x^{\alpha}).
\end{eqnarray*}

\end{ex}
In this example, we have connected two important functions the
exponential function in the ordinary space and
$\alpha$-exponential function in the $\alpha$-fractional space.

 \begin{center}\textbf{CORRESPONDENCE}\end{center}

\begin{tabular}{|c|c|c|}
  \hline
  \mbox{ordinary space} &  & \mbox{$\alpha$-level space}  \\
  \hline
  $1$ &  & $\frac{x^{\alpha-1}}{\Gamma(\alpha)}$ \\
\hline $x$&&$\frac{x^{2\alpha-1}}{\Gamma(2\alpha)}$\\
  \hline
  $ e^{-x}$ &  & $x^{\alpha-1}E_{(\alpha,\alpha)} (- x^{\alpha}) $ \\
  \hline
$ e^{-ax}e^{-bx},~ e^{-(a+b)x}$ &  & $x^{\alpha-1}E_{(\alpha,\alpha)} (- (a+b)x^{\alpha}) $ \\
\hline
   $\frac{x^{\eta-1}e^{-\frac{x}{\delta}}}{\delta^{\eta}\Gamma(\eta)}$&  & $\frac{x^{\alpha\eta-1}}{\delta^{\eta}}E^{\eta}_{(\alpha,\alpha\eta)}
(-\delta^{-1} x^{\alpha})$ \\
\hline
  $_{0}F_{1}(;a;-t)$ &  & $x^{\alpha-1}\sum^{\infty}_{k=0}\frac{ (-1)^{k} x^{\alpha
k } }{(a)_{k} \Gamma(\alpha k+\alpha)},~~a > 1$ \\
\hline

 $_{1}F_{0}(a;;-t)$ &  & $x^{\alpha-1}\sum^{\infty}_{k=0}\frac{ (a)_{k}(-1)^{k} x^{\alpha
k } }{ \Gamma(\alpha k+\alpha)},~~a > 1$ \\
\hline
%
$_{1}F_{1}(a;b;-t)$ &  & $x^{\alpha-1}\sum^{\infty}_{k=0}\frac{
(a)_{k}(-1)^{k} x^{\alpha
k } }{(b)_{k} \Gamma(\alpha k+\alpha)},~~a > 1, b > 1$ \\
\hline
 $_{2}F_{1}(a,b;c;-t)$ &  &
$x^{\alpha-1}\sum^{\infty}_{k=0}\frac{ (a)_{k}(b)_{k}(-1)^{k}
x^{\alpha
k } }{(c)_{k} \Gamma(\alpha k+\alpha)},~~a > 1, b>1, c>1$ \\
  \hline
\end{tabular}

\section{Conclusion}\label{s5}

 It seems that the Levy structure is important in the area
of Fractional Calculus in the following sense,

\begin{itemize}
    \item fractional equations
naturally possess the Levy structure such as \eqref{hhh1}
    \item functions with Levy structure matches fittingly with
    functions with Levy structure, for example the tables shown in section 3
\end{itemize}

From section 3 and 4, it is urged that the theory of fractional
calculus be developed  according to the origin of functions where
the origin means the level of functions such as $\alpha$-level.

For this, see \cite{ki1, mh3}.\\

 Furthermore, if we use the asymptotic behavior of the gamma
function , $x^{\alpha-1}E_{(\alpha,\alpha)}(-ax^{\alpha})$ can be
approximately used in place of $\sum^{\infty}_{k=0}\frac{(a)_{k}
(-1)^{k} x^{\alpha k +\alpha-1} }{ \Gamma(\alpha k+\alpha)} $,
namely, the coefficient of $x$ can be generated by using the
limiting process \cite{mh2}.
\\\\
\textbf{\large Acknowledgement}\\ The author would like to name
the transform "Mathai transform" in section 4 after the Emeritus
Professor Arakaparampil M. Mathai to honour him and his
contributions in these areas. The author would like to thank to
the Department of Science and Technology, Government of India, New
Delhi, for the financial assistance under Project No.
SR/S4/MS:287/05 and the Centre for Mathematical Sciences, Pala
campus, for providing all facilities.


\begin{thebibliography}{7}

\bibitem{ga}
 Gamelin T. W.,
  Complex Analysis, Springer, 2004.

\bibitem{ki1}
 Kilbas A. A., Srivastava H. M.,  Trujillo J. J.,
  Theory and Applications of Fractional Differential Equations,
  Elsevier, 2006.

\bibitem{ma3}
Mathai A.M., A Handbook of Generalized Special Functions for
Statistical and Physical Sciences, Oxford University Press,
Oxford, 1993.

\bibitem{ma2}
 Mathai A.M., A pathway to matrix-variate gamma and Gaussian
densities, Linear Algebra and Its Applications, 396(2005),
317-328.

 \bibitem{ma4}
 Mathai A.M., Basic Probability and Statistics, part 1, Probability and
 random variables, module 6, Centre for mathematical Sciences and printed
 at Mathematical Sciences Press, CMS Pala Campus, India, 2010.

\bibitem{ma1}
 Mathai A.M., Some properties of Mittag-Leffler functions and
matrix-variate analogues: a statistical perspective, Preprint.




\bibitem{mh1}
 Mathai A.M. and  Haubold Hans J., A general overview of pathway
model, Tsallis statistics and generalizations, Preprint.

\bibitem{mh2}
Mathai A.M. and  Haubold Hans J., Mittag-Leffler functions to
Pathway Model to Tsallis statistics, Preprint.

\bibitem{mh3}
Mathai A.M. and  Haubold Hans J., Special Functions for Applied
Scientists, Springer, New York, 2008.


\bibitem{mhs2}
Mathai A.M., Saxena R.K. and  Haubold Hans J., The H-function:
Theory and Applications, Springer, New York, 2010.

\end{thebibliography}
\end{document}